\documentclass{article}

\usepackage{latexsym}
\usepackage{amsfonts,amssymb} 
\usepackage{graphicx}

\def\cal{\mathcal}

\newcommand{\field}[1]{\mathbb{#1}}
\newcommand{\A}{\field{A}}
\newcommand{\C}{\field{C}}

\newcommand{\N}{\field{N}}

\newcommand{\R}{\field{R}}

\newcommand{\Z}{\field{Z}}

\newcommand{\pgcd}{{\rm pgcd}}
\newcommand{\kk}{{\bf k}}

\newtheorem{theorem}{Theorem}[section]
\newtheorem{proposition}{Proposition}[section]
\newtheorem{lemma}{Lemma}[section]
\newtheorem{corollary}{Corollary}[section]
\newtheorem{definition}{Definition}[section]
\newtheorem{remark}{Remark}[section]
\newtheorem{conjecture}{Conjecture}[section]

\begin{document}

\makeatletter	   
\makeatother     

\title{Curves Defined by Chebyshev Polynomials}         
\author{Gene Freudenburg and Jenna Freudenburg}
\date{\today} 

\maketitle

\pagestyle{plain}

\begin{abstract} Working over a field $\kk$ of characteristic zero, this paper studies line embeddings of the form $\phi = (T_i,T_j,T_k):\A^1\to\A^3$, where $T_n$ denotes the degree $n$ Chebyshev polynomial of the first kind. In {\it Section 4}, it is shown that (1) $\phi$ is an embedding if and only if the pairwise greatest common divisor of $i,j,k$ is 1, and (2) for a fixed pair $i,j$ of relatively prime positive integers, the embeddings of the form $(T_i,T_j,T_k)$ represent a finite number of algebraic equivalence classes. {\it Section 2} gives an algebraic definition of the Chebyshev polynomials, where their basic identities are established, and {\it Section 3} studies the plane curves $(T_i,T_j)$. {\it Section 5} establishes the Parity Property for Nodal Curves, and uses this to parametrize the family of alternating $(i,j)$-knots over the real numbers. 
\end{abstract}

\section{Introduction} 

Every knot $K\subset S^3$ can be parametrized by polynomials on the open set $\R^3\subset S^3$, where 
$\R^3 = S^3 -\{ P\}$ for some point $P\in K$. This was shown by Shastri \cite{Shastri.92} using the Weierstrass Approximation Theorem. However, finding a parametrization for a specific knot $K$ is difficult. The standard method has been to first find a parametrization $(f(t),g(t))$ of a regular plane projection of $K-\{ P\}$, and to then find a third polynomial $h(t)$ giving the correct over- and under-crossings at the nodes of the projection. Consequently, the list of knots for which explicit polynomial parametrizations have been found to date is finite.

\begin{definition} {\rm A set of knots $\cal K$ is {\it topologically infinite}, or {\it t-infinite}, if and only if $\cal K$ contains a sequence $K_i$, $i\ge 0$, such that $K_i$ and $K_j$ are distinct knot types when $i\ne j$. Otherwise, $\cal K$ is {\it t-finite}, i.e., elements of $\cal K$ represent a finite number of knot types. }
\end{definition} 
The results presented in this paper include the following.

\begin{itemize}
\item[1.] 
Let $T_n$ denote the degree-$n$ Chebyshev polynomial of the first kind over a field $\kk$ of characteristic zero. Given positive integers $i,j,k$, define their {\it pairwise greatest common divisor} to be 
\[
\pgcd (i,j,k)= \max\{ \gcd (i,j), \gcd (i,k), \gcd (j,k)\}\,\, .
\]
{\it Proposition~\ref{embedding}} asserts that $\phi =(T_i,T_j,T_k):\A^1\to\A^3$ is an algebraic line embedding if and only if $\pgcd (i,j,k)=1$. 
\item[2.] {\it Proposition~\ref{reduced}} shows that, for a fixed pair of positive integers $(i,j)$ with 
$\gcd (i,j)=1$, the infinite set
\[
{\cal K}_{(i,j)}:=\{ (T_i,T_j,T_k)\,\,\vert\,\, k\ge 1, \pgcd (i,j,k)=1\}
\]
represents a finite number of algebraic equivalence classes. This implies that ${\cal K}_{(i,j)}$ is t-finite over the field $\R$ of real numbers. 
\item[3.] {\it Proposition~\ref{t-infinite}} gives an explicit polynomial parametrization for each member of the t-infinite family of alternating $(i,j)$-knots, using Chebyshev polynomials of both the first and second kind. 
\end{itemize} 
The paper is organized in the following way. 
\begin{itemize}
\item Section 2: The Chebyshev Polynomials $T_n$
\item Section 3: The Curves $(T_i,T_j)$
\item Section 4: Line Embeddings $(T_i,T_j,T_k)$
\item Section 5: Nodal Curves and Crossing Sequences
\item Section 6: Open Problems
\end{itemize}

The literature on polynomial line embeddings in 3-space is scant. Following the appearance of the Epimorphism Theorem in the mid-1970s, which characterizes embeddings of the line in a plane, Abhyankar considered the family of embeddings of the form $\theta (\ell ,m,n)=(t+t^{\ell},t^m,t^n)$ for positive integers $\ell ,m,n$ \cite{Abhyankar.77}. Abhyankar conjectured that these embeddings are not rectifiable when none of the three integers belongs to the semigroup generated by the other two . However, Craighero showed that $\theta (5,4,3)$ and $\theta (6,5,4)$ can be rectified \cite{Craighero.86,Craighero.88}; see \cite{Freudenburg.06}, 4.3.1, for another proof. Subsequently, Bhatwadekar and Roy gave families of integers $\ell ,m,n$ for which $\theta (\ell ,m,n)$ can be rectified \cite{Bhatwadekar.Roy.91}. The general case remains open, including $\theta (7,6,5)$.  

The lines defined by $\theta (\ell ,m,n)$ are topologically trivial over the field of real numbers. In a well-known paper \cite{Vassiliev.90}, Vassiliev observed that any ``non-compact'' knot can be represented by a polynomial embedding. Shastri \cite{Shastri.92} gave a proof of this fact and provided specific polynomial embeddings for the trefoil and figure-8 knots. Thus, the Gordon-Lueke Theorem from topology implies that there exist algebraic embeddings of $\R^1$ in $\R^3$ which are {\it not} rectifiable. It is an important open question whether algebraic embeddings of $\C^1$ in $\C^3$ can always be rectified by algebraic automorphisms; Kaliman \cite{Kaliman.92} showed that they can be rectified by holomorphic automorphisms. 

Most work on the subject beyond that described above has been aimed at finding parametrizations for specific knots, or determining minimal degrees for such parametrizations \cite{Koseleff.Pecker.ppta, Koseleff.Pecker.pptb, Mishra.99, Mishra.00, Ranjan.Shukla.96}. In their recent preprint \cite{Durfee.Oshea.ppt}, Durfee and O'Shea survey results on polynomial knots, including parametrizations for several specific knots. 

The results from {\it Sect. 4} of this paper were reported previously in the two research reports \cite{J.Freudenburg.08a, J.Freudenburg.08b} of the second author. 

\paragraph{Preliminaries.}
The ground field $\kk$ is any field of characteristic zero, though the field of real numbers $\kk =\R$ is used in discussing knots. The polynomial ring in $n$ variables $x_1,...,x_n$ is denoted by $\kk^{[n]}=\kk [x_1,...,x_n]$. The subalgebra generated by $f_1,...,f_m\in\kk [x_1,...,x_n]$ is denoted by $\kk [f_1,...,f_m]$. 
Affine $n$-space over $\kk$ is $\A^n$. 

Let $\phi :\A^1\to\A^n$ be a morphism of $\kk$-varieties, and let $\phi^*:\kk [x_1,...,x_n]\to\kk [t]$ be the corresponding map of \kk -algebras. The following statements are equivalent (see \cite{Essen.00}, Cor. B.2.6, or \cite{Durfee.Oshea.ppt}, Lemma 4).
\begin{itemize}
  \item[1.] $\phi$ is an embedding.
  \item[2.] $\phi$ is injective, and $\phi '(t)\ne 0$ for all $t\in\A^1$.
  \item[3.] $\phi *$ is surjective.
\end{itemize}  

$GA_n(k)$ is the group of algebraic automorphisms of $\kk [x_1,...,x_n]$, or equivalently, of $\A^n$. Two embeddings
$\phi ,\psi :\A^1\to\A^n$ are said to be {\it algebraically equivalent} if and only if there exists $\alpha\in GA_n(\kk )$ with $\alpha\phi=\psi$. An embedding $\phi :\A^1\to\A^n$ is {\it rectifiable} if and only if $\phi$ is algebraically equivalent to the standard embedding defined by $t\to (t,0,0)$. The famous Epimorphism Theorem of Abhyankar and Moh \cite{Abhyankar.Moh.75} and Suzuki \cite{Suzuki.74} asserts that every algebraic embedding of $\A^1$ in $\A^2$ is algebraically equivalent to the standard embedding. We use the following equivalent version.
\begin{theorem} {\rm (Epimorphism Theorem)} Assume $\kk$ is a field of characteristic zero. If $f(t),g(t)\in\kk [t]$ are such that $t\in\kk [f(t),g(t)]$, then either $\deg f\vert\deg g$ or $\deg g\vert\deg f$.
\end{theorem}
Many other proofs of this result have subsequently appeared. Of particular note is the paper of Rudolph \cite{Rudolph.82}, which gives an elegant proof of the Epimorphism Theorem using knot theory. 

Shastri \cite{Shastri.92} gives the embedding $\phi :\R^1\to\R^3$ defined by $\phi (t)=(t^3-3t,t^4-4t^2,t^5-10t)$, and shows that $\phi (\R^1)$ is the trefoil knot. In particular, if $F\in k[X,Y,Z]$ is $F=YZ-X^3-5XY+2Z-7X$, then $\phi^*(F)=t$. While $\phi$ itself is already quite simple, we see that $F$ can be simplified:
$F=(Y+2)(Z-5X)+3X-X^3$. Thus, if $\alpha\in GA_3(\R )$ is defined by $\alpha =(X,Y-2,Z+5X)$, then 
\[
\alpha\phi = (t^3-3t, t^4-4t^2+2, t^5-5t^3+5t) \quad {\rm and} \quad \alpha (F)=YZ+3X-X^3\,\, .
\]
We recognize that the defining polynomials for the trefoil in this re-parametrization are the (monic) Chebyshev polynomials of degree 3,4, and 5. It is therefore natural to ask if other combinations of Chebyshev polynomials define embeddings, and if so, what knots they parametrize over the field of real numbers. 

In general, note that two kinds of equivalence of algebraic lines in $\A^3$ are considered in this paper. First, the notion of algebraic equivalence is defined above, and is valid for any ground field. In addition, when the ground field is $\kk =\R$, we consider topological equivalence, by which we mean the existence of a homeomorphism from $\R^3$ to itself which carries one embedded line to the other. Thus, a knot is not distinguished from its mirror image, as with isotopic equivalence. Clearly, in the case of algebraic lines in $\R^3$, algebraic equivalence implies topological equivalence. 

\paragraph{Note Added in Proof.} The paper {\it Chebyshev Knots} \cite{Koseleff.Pecker.pptc} of Koseleff and Pecker was posted on the arXiv preprint server in December 2008. Their paper also studies knots parametrized by Chebyshev polynomials, and several of their results parallel those found in our paper. In some cases, their more topological approach enables them to give proofs of results we had conjectured, most notably: (1) for $n\ge 1$, $(T_3,T_{3n+1},T_{3n+2})$ parametrizes the $(2,2n+1)$ torus knot, and (2) for relatively prime $i$ and $j$, $(T_i,T_j, T_{ij-i-j})$ parametrizes the alternating $(i,j)$-knot (which they call an alternate harmonic knot). On the other hand, the focus of our paper is on algebraic aspects, for example, working over a field of characteristic zero when possible, with empasis on algebraic equivalence/non-equivalence of embedded lines. 

\section{The Chebyshev Polynomials $T_n$}
\subsection{The Polynomials $f_n$ and $g_n$}
Given $P\in\Z^{[2]}$, there exist $Q,R\in\Z^{[2]}$ such that
$P(x,iy)=Q(x,y^2)+iyR(x,y^2)$, where $i^2=-1$. In particular, for any integer $n\ge 0$, let $f_n\in\Z[x,y^2]$ and $g_n\in y\Z [x,y^2]$ be such that 
\[
(x+iy)^n = f_n+ig_n \,\, .
\]
Note that $f_n$ and $g_n$ are homogeneous of degree $n$. Note also that replacing $y$ by $(-y)$ gives 
\[
(x-iy)^n=f_n-ig_n\,\, .
\]
\begin{lemma}\label{basic} Assume $m,n\in\Z$ and $m\ge n\ge 0$. 
\begin{itemize}
  \item[{\bf (a)}] $f_{m+n}=f_mf_n-g_mg_n$ \quad and\quad  $g_{m+n}=g_mf_n+f_mg_n$
  \item[{\bf (b)}] $(x^2+y^2)^nf_{m-n}=f_mf_n+g_mg_n$ \quad and\quad  $(x^2+y^2)^ng_{m-n}=g_mf_n-f_mg_n$
  \item[{\bf (c)}] $(x^2+y^2)^n=f_n^2+g_n^2$
  \item[{\bf (d)}] $f_n(f_m,g_m)=f_m(f_n,g_n)=f_{mn}$ \quad and\quad  $g_n(f_m,g_m)=g_m(f_n,g_n)=g_{mn}$
\end{itemize}  
\end{lemma}

\noindent {\it Proof.} Statement (a) follows immediately by comparing $(x+iy)^m(x+iy)^n$ with $(x+iy)^{m+n}$. For (b), we have:
\begin{eqnarray*}
(x^2+y^2)^n(f_{m-n}+ig_{m-n}) &=& (x^2+y^2)^n(x+iy)^{m-n}\\
&=& (x+iy)^n(x-iy)^n(x+iy)^{m-n}\\
&=& (x+iy)^m(x-iy)^n\\
&=& (f_m+ig_m)(f_n-ig_n)\\
&=& (f_mf_n+g_mg_n)+i(g_mf_n-f_mg_n)\,\, .
\end{eqnarray*}
Part (c) is a special case of (b) when $m=n$, since $f_0=1$. Finally, part (d) follows immediately by comparing $\bigl[ (x+iy)^m\bigr]^n$ with $(x+iy)^{mn}$. $\Box$

\subsection{Chebyshev Polynomials over $\Z$}

\begin{definition} {\rm Define $\pi :\Z [x,y^2]\to\Z[t]$ by $\pi (x)=t$ and $\pi (y^2)=1-t^2$. Define $T_n,U_n\in\Z [t]$ by
\[
T_n=\pi (f_n)\quad  (n\ge 0)\quad  {\rm and}\quad  U_n=\pi (g_{n+1}/y)\quad (n\ge -1)\,\,  .
\]
Then $T_n$ and $U_n$ are the} Chebyshev polynomials {\rm of the first and second kind, respectively.}
\end{definition}
In particular, note that $U_{-1}=0$, $T_0=U_0=1$, $T_1=t$, and $U_1=2t$. As with trigonometric functions, there are countless identities involving Chebyshev polynomials. Those needed in this paper are given in the following lemma.
\begin{lemma}\label{identities} Given $m,n\in\Z$ with $m\ge n\ge 0$:
\begin{itemize}
\item[{\bf (a)}] $T_{m+n}=T_mT_n-(1-t^2)U_{m-1}U_{n-1}$
\item[{\bf (b)}] $T_{m-n}=T_mT_n+(1-t^2)U_{m-1}U_{n-1}$
\item[{\bf (c)}] $2T_mT_n=T_{m+n}+T_{m-n}$
\item[{\bf (d)}] $U_{m+n}=U_mT_n+T_{m+1}U_{n-1}$
\item[{\bf (e)}] $U_{m-n}=U_mT_n-T_{m+1}U_{n-1}$
\item[{\bf (f)}] $2U_mT_n=U_{m+n}+U_{m-n}$
\item[{\bf (g)}] $T_m(T_n)=T_n(T_m)=T_{mn}$
\item[{\bf (h)}] $U_{mn-1}(t)=U_{n-1}(T_m)U_{m-1}=U_{m-1}(T_n)U_{n-1}$
\item[{\bf (i)}] $T_n'=nU_{n-1}$ 
\end{itemize}
\end{lemma}

\medskip

\noindent {\it Proof.} {\it Lemma~\ref{basic}(a)} implies 
\[
T_{m+n}=\pi (f_mf_n-g_mg_n)=\pi (f_m)\pi (f_n)-\pi (y^2)\pi (g_m/y)\pi (g_n/y)=T_mT_n-(1-t^2)U_{m-1}U_{n-1} \,\, .
\]
In addition, since $\pi (x^2+y^2)=1$, {\it Lemma~\ref{basic}(b)} implies
\[
T_{m-n}=\pi (f_mf_n+g_mg_n)=\pi (f_m)\pi (f_n)+\pi (y^2)\pi (g_m/y)\pi (g_n/y)=T_mT_n+(1-t^2)U_{m-1}U_{n-1} \,\, .
\]
Parts (a) and (b) are thus proved, and part (c) follows by adding these two equalities. 

Similarly, {\it Lemma~\ref{basic}(a)} implies
\[
U_{m+n}=\pi ((g_{m+1}/y)f_n+f_{m+1}(g_n/y)) =\pi (g_{m+1}/y)\pi (f_n)+\pi (f_{m+1})\pi (g_n/y)=U_mT_n+T_{m+1}U_{n-1} \,\, .
\]
In addition, since $\pi (x^2+y^2)=1$, {\it Lemma~\ref{basic}(b)} implies 
\[
U_{m-n}=\pi ((g_{m+1}/y)f_n-f_{m+1}(g_n/y))=\pi (g_{m+1}/y)\pi (f_n)-\pi (f_{m+1})\pi (g_n/y)=U_mT_n-T_{m+1}U_{n-1} \,\, .
\]
Parts (d) and (e) are thus proved, and part (f) follows by adding these two equalities. 

In order to prove (g), for each $n\ge 0$ let $F_n,G_n\in\Z^{[2]}$ be such that $f_n(x,y)=F_n(x,y^2)$ and $g_n(x,y)=yG_n(x,y^2)$, 
noting that $T_n(t)=F_n(t,1-t^2)$ and $U_n(t)=G_{n+1}(t,1-t^2)$. Since $g_n^2\in\Z [x,y^2]$ for all $n\ge 0$, {\it Lemma~\ref{basic}(c)} implies
\[
1=\pi (f_n^2+g_n^2)=T_n^2+\pi (g_n^2)\,\, .
\]
Therefore, by {\it Lemma~\ref{basic}(d)}, it follows that
\[
T_{mn}(t)=\pi \big( f_n (f_m,g_m)\bigr) =\pi \big( F_n(f_m,g_m^2)\bigl)=F_n \bigl( T_m(t),1-T_m^2(t))\bigr) = T_n (T_m(t))\,\, .
\]
By symmetry, $T_{mn}(t)=T_m (T_n(t))$ as well, so (g) is proved.

For (h), it follows from {\it Lemma~\ref{basic}(d)} that
\begin{eqnarray*}
U_{mn-1}(t) &=& \pi\bigl( g_{mn}/y\bigr) \\
&=& \pi\bigl( g_n(f_m,g_m)/y\bigr) \\
&=& \pi\bigl( g_mG_n(f_m,g_m^2)/y\bigr) \\
&=& \pi\bigl( G_n(f_m,g_m^2)\bigr)\cdot\pi\bigl( g_{m}/y\bigr) \\
&=& G_n\bigl( T_m(t),1-T_m(t)^2\bigr)\cdot U_{m-1}(t)\\
&=& U_{n-1}\bigl( T_m(t)\bigr) U_{m-1}(t)\,\, .
\end{eqnarray*}
The remaining equality of (h) follows by symmetry.

Finally, to prove part (i) we need the following three identities: 
\[
T_{n+1}=2tT_n-T_{n-1}\,\, ,\quad U_n=2tU_{n-1}-U_{n-2} \quad {\rm and} \quad U_n=2T_n+U_{n-2}\quad (n\ge 1)\,\, .
\]
The first of these is implied by part (a) when $m=1$, and the second follows from (d) when $n=1$. 
For the third, set $m=0$ in part (d) to obtain $U_n=T_n+tU_{n-1}$. Thus,
\[
2T_n=2U_n-2tU_{n-1}=2U_n-(U_n+U_{n-2})=U_n-U_{n-2}\,\, .
\]
We now prove (i) by induction, the cases $n=0,1$ being clear. Assume $n\ge 1$ and $T_m'=mU_{m-1}$ whenever $0\le m\le n$. Then
\begin{eqnarray*}
T_{n+1}' &=& 2t\cdot T_n'+2T_n-T_{n-1}'\\
&=& 2t\cdot nU_{n-1} + 2T_n-(n-1)U_{n-2}\\
&=& n(2tU_{n-1}-U_{n-2})+(2T_n+U_{n-2})\\
&=& nU_n+U_n \\
&=& (n+1)U_n
\end{eqnarray*}
This completes the proof of the lemma. $\Box$

\medskip

Note that as a special case, identities (a) and (d) above give the familiar recursion formulae:
\[
T_0=1\quad ,\quad T_1=t\quad {\rm and}\quad T_{n+1}=2tT_n-T_{n-1} \quad (n\ge 1)
\]
and
\[
U_{-1}=0\quad ,\quad U_0=1\quad {\rm and}\quad U_{n+1}=2tU_n-U_{n-1} \quad (n\ge 0)
\]
From this, it is clear that $\deg T_n=n$, and that the leading coefficient of $T_n$ is $2^{n-1}$. In addition, $T_n$ is an even function for even $n$, and an odd function for odd $n$. By induction, these also show that $T_n(1)=1$ and $U_n(1)=n+1$. 

The Chebyshev polynomials $T_n(t)$ share many properties with the monomials $t^n$. They form a basis for $\kk [t]$ as a $\kk$-vector space. Moreover, up to a certain mild equivalence, these are the only families of univariate polynomials which commute with each other by composition: 
\[
t^m\circ t^n=t^n\circ t^m=t^{mn} \quad {\rm and}\quad T_m\circ T_n=T_n\circ T_m=T_{mn}\,\, .
\]
See \cite{Schinzel.82}.

\subsection{Chebyshev Polynomials over $\R$}

\begin{lemma}\label{roots} Let $n$ be a non-negative integer. 
\begin{itemize}
  \item[{\bf (a)}] For all $\theta\in\R$,
\[
T_n(\cos\theta )=\cos (n\theta ) \quad {\rm and} \quad U_n(\cos\theta )=\frac{\sin ((n+1)\theta )}{\sin\theta}\,\, .
\]
  \item[{\bf (b)}] The roots of $T_n$ are $\cos\bigl(\frac{2k-1}{2n}\pi\bigr)$, $k=1,...,n$, and each of these is a simple root.
  \item[{\bf (c)}] The roots of $U_n$ are $\cos\bigl(\frac{k}{n+1}\pi\bigr)$, $k=1,...,n$, and each of these is a simple root.
  \item[{\bf (d)}] Given $t\in \R$, $-1\le t\le 1$ if and only if $-1\le T_n(t)\le 1$.
\end{itemize}  
\end{lemma}

\noindent {\it Proof.} Let $F_n,G_n\in\Z^{[2]}$ be defined as in the preceding proof, i.e., $f_n(x,y)=F_n(x,y^2)$ and $g_n(x,y)=yG_n(x,y^2)$. Then
\[
T_n(t)=F_n(t,1-t^2) \quad {\rm and}\quad U_n(t)= G_{n+1}(t,1-t^2)\,\, .
\]
It follows that, for any real number $\theta$, 
\begin{eqnarray*}
\cos (n\theta )+i\sin (n\theta ) &=& (\cos\theta + i\sin\theta)^n \\
&=& f_n(\cos\theta ,\sin\theta ) + i g_n(\cos\theta ,\sin\theta )\\
&=& F_n(\cos\theta ,\sin^2\theta )+i\sin\theta G_n(\cos\theta ,\sin^2\theta )\\
&=& F_n(\cos\theta ,1-\cos^2\theta )+i\sin\theta G_n(\cos\theta ,1-\cos^2\theta )\\
&=& T_n(\cos\theta) + i\sin\theta U_{n-1}(\cos\theta)\,\, .
\end{eqnarray*} 
This proves (a). Next, note that 
\[
\Bigl\lbrace\cos \Bigl(\frac{2k-1}{2n}\pi\Bigr)\,\,\vert\,\, k=1,...,n\Bigr\rbrace
\]
is a set of $n$ distinct real numbers, and from part (a), each number in this set is a root of $T_n(t)$. Since $\deg T_n=n$, part (b) follows. A similar argument is used to show (c). 

For (d), let $t\in\R$ be given, $-1\le t\le 1$, and choose $\theta\in\R$ such that $t=\cos\theta$. Then by part (a), it follows that $-1\le T_n(t)\le 1$. Conversely, note that $T_n$ can have no critical points outside the interval $[-1,1]$, since $T_n'=nU_{n-1}$ and the roots of $U_{n-1}$ lie in the interval $[-1,1]$. Therefore, since $U_n(1)=n+1>0$, $T_n$ is increasing for $t>1$, and since $T_n(1)=1$, this implies $T_n(t)>1$ for $t>1$. Finally, 
since $\vert T_n(-t)\vert = \vert T_n(t)\vert$, it follows that $\vert T_n(t)\vert > 1$ when $\vert t\vert >1$. $\Box$

\section{The Curves $(T_i,T_j)$}

Given positive integers $i$ and $j$ with $\gcd{(i,j)}=1$, define $F\in k[x,y]$ by $E:=T_i(y)-T_j(x)$. Consider the degree function on $\kk [x,y]$ defined by 
$\deg x= i$ and $\deg y=j$. Given $f\in\kk [x,y]$, let $\bar{f}$ denote the highest homogeneous summand of $f$ relative to this grading. Clearly, if $F=GH$ for $F,G\in\kk [x,y,]$ of positive degree, then
$\bar{F}=\bar{G}\bar{H}$, where $\bar{G}$ and $\bar{H}$ are of positive degree. However, $\bar{F}=2^{i-1}y^i-2^{j-1}x^j$, and it is well known that this polynomial is irreducible when $i$ and $j$ are coprime. Therefore, $F$ is irreducible. 

Let $C$ be the irreducible curve defined by $F=T_i(y)-T_j(x)=0$. 
Since $T_i\circ T_j=T_j\circ T_i$, it follows that $C$ is parametrized by $\alpha =(T_i,T_j):\A^1\to\A^2$. 
Note that by {\it Lemma~\ref{identities}(i)}, $\alpha^{\prime}(t)=(iU_{i-1},jU_{j-1})$. By {\it Lemma~\ref{roots}(c)}, it follows that $\alpha^{\prime}(t)\ne 0$ for all $t\in k$. 

\begin{proposition}\label{curves} Assume $i$ and $j$ are positive integers with $\gcd (i,j)=1$, and let $C$ be the plane curve defined by $T_i(y)-T_j(x)=0$. 
\begin{itemize}
\item[{\bf (a)}] The singularities of $C$ consist of $\frac{1}{2}(i-1)(j-1)$ nodes.
\item[{\bf (b)}] The nodes of $C$ are precisely points of the set
\[
S= \Bigl\lbrace \Bigl( \cos \Bigl(\frac{\lambda\pi}{j}\Bigr)\, ,\, 
\cos\Bigl(\frac{\mu\pi}{i}\Bigr)\Bigr)\,\, \Bigl| \,\, 1\le\lambda\le j-1\, ,\,\, 1\le\mu\le i-1\, ,\,\, \lambda\equiv\mu\, ({\rm mod}\, 2)
\Bigr\rbrace \,\, .
\]
\item[{\bf (c)}] Given the node 
$\displaystyle P=\Bigl( \cos\Bigl(\frac{\lambda\pi}{j}\Bigr)\, ,\, \cos\Bigl(\frac{\mu\pi}{i}\Bigr)\Bigr)$, 
\[
\alpha^{-1}(P)=\Bigl\lbrace \cos \Bigl(\frac{k_1\pi}{ij}\Bigr)\,\, ,\,\,
\cos\Bigl(\frac{k_2\pi}{ij}\Bigr)\Bigr\rbrace\,\, ,
\]
where $k_1 = \lambda iu + \mu jv$, $k_2 = \lambda iu - \mu jv$, and $u,v\in\Z$ satisfy $iu+jv=1$.

\item[{\bf (d)}] 
\[
\alpha^{-1}(S)=\Bigl\lbrace\,\,\cos \Bigl(\frac{k\pi}{ij}\Bigr)\,\,\Bigl|
\,\, 1\le k<ij\, ,\, k\not\in i\Z\, ,\, k\not\in j\Z\,\,\Bigr\rbrace\,\, .
\]
\end{itemize}
\end{proposition}

\noindent {\it Proof.} 
If $F(x,y)=T_i(y)-T_j(x)$, then the singular points of $C$ are defined by the system
\[
\left\{ \begin{array}{llll}
F(x,y)&=&T_i(y)-T_j(x)&=0\\
F_x(x,y)&=&-T_j'(x)&=0\\
F_y(x,y)&=& T_i'(y)&=0\,\, .
\end{array}\right.
\]
By {\it Lemma~\ref{identities}(k)} and {\it Lemma~\ref{roots}(c)}, the roots of $T_j'(x)$ are $\cos (\frac{k\pi}{j})$ $(k=1,...,j-1)$, and the roots of $T_i'(x)$ are $\cos (\frac{k\pi}{i})$, $(k=1,...,i-1)$. It follows that the set of singular points of $C$ is $S\cap C$. 

Let $P=(a,b)$ be a point of $S$, where 
\[
a=\cos \Bigl(\frac{\lambda\pi}{j}\Bigr) \quad  {\rm and}\quad  b=\cos \Bigl(\frac{\mu\pi}{i}\Bigr)
\]
as above. Then
\begin{eqnarray*}
F(P)&=&T_i\bigl( \cos ({\textstyle\frac{\mu\pi}{i}})\bigr)-T_j\bigl( \cos ({\textstyle\frac{\lambda\pi}{j}})\bigr)\\
    &=&\cos \bigl( i{\textstyle \frac{\mu\pi}{i}}\bigr)-\cos \bigl( j{\textstyle \frac{\lambda\pi}{j}}\bigr)\\
    &=&\cos (\mu\pi )-\cos (\lambda\pi )\\
    &=&0\,\, ,
\end{eqnarray*}
since $\lambda\equiv\mu\, ({\rm mod}\, 2)$. Therefore, $S\subset C$, and $S$ is precisely the set of singular points of $C$. 

It remains to show that $P$ is an ordinary double point of $C$. 
For $X=x-a$ and $Y=y-b$, let $G\in\kk^{[2]}$ be such that 
$G(X,Y)=F(x,y)$, and write $G(X,Y)=\sum_{d\ge 0} G_d(X,Y)$, where $G_d$ is homogeneous of degree $d$ in $X$ and $Y$.
Then
\[
G(X,Y)=T_i(y)-T_j(x)=T_i(Y+b)-T_j(X+a)=\sum_{m=0}^i\frac{T_i^{(m)}(b)}{m!}Y^m - \sum_{n=0}^j\frac{T_j^{(n)}(a)}{n!}X^n\,\, .
\]
In particular, 
\[
G_2={\textstyle\frac{1}{2}}\bigl( T_i''(b)Y^2-T_j''(a)X^2\bigr)\,\, .
\]
Note that $F_y(P)=T_i'(b)=0$ and $F_x(P)=-T_j'(a)=0$.
Since the roots of $T_i'(t)=iU_{i-1}(t)$ and $T_j'(t)=jU_{j-1}(t)$ are simple, it follows that 
$T_i''(b)\ne 0$ and $T_j''(a)\ne 0$. Therefore, $G_2$ factors as the product of two distinct linear forms. It follows that $P$ is a node of $C$. So parts (a)  and (b) are proved.

In order to prove part (c), let $\alpha=(T_i,T_j)$, and set $t_m=\cos\bigl( \frac{k_m\pi}{ij}\bigr)$ for $k_m$ defined above ($m=1,2$). If $t_1=t_2$, then
\[ 
\frac{k_1\pi}{ij}=\pm\frac{k_2\pi}{ij}+2N\pi\quad {\rm for\,\, some}\,\, N\in\Z \quad \Rightarrow\quad 
\lambda iu+\mu jv=\pm(\lambda iu-\mu jv)+2ijN\,\, .
\]
Thus, either $2\mu jv=2ijN$ (positive case) or $2\lambda iu=2ijN$ (negative case). In the first case,  the fact that $iu+jv=1$ implies $\mu\equiv 0\, ({\rm mod}\, i)$, which is impossible since $1\le\mu\le i-1$. Likewise, the second case yields 
$\lambda\equiv 0\, ({\rm mod}\, j)$, which is impossible since $1\le\lambda\le j-1$. Therefore, $t_1\ne t_2$. 

On the other hand, $\alpha (t_1)=\alpha (t_2)$, since for $m=1,2$, 
\begin{eqnarray*}
\alpha (t_m)&=&\Bigl( T_i\bigl(\cos {\textstyle\frac{k_m\pi}{ij}}\bigr)\, ,\, T_j\bigl(\cos {\textstyle\frac{k_m\pi}{ij}}\bigr)\Bigr)\\
            &=&\Bigl(\cos\bigl(i{\textstyle\frac{k_m\pi}{ij}}\bigr)\, ,\, \cos\bigl(j{\textstyle\frac{k_m\pi}{ij}}\bigr)\Bigr)\\
            &=&\Bigl(\cos\bigl({\textstyle\frac{(\lambda iu\pm\mu jv)\pi}{j}}\bigr)\, ,\, 
                     \cos\bigl({\textstyle\frac{(\lambda iu\pm\mu jv)\pi}{i}}\bigr)\Bigr)\\
            &=&\Bigl(\cos\bigl({\textstyle\frac{(\lambda (1-jv)\pm\mu jv)\pi}{j}}\bigr)\, ,\, 
                     \cos\bigl({\textstyle\frac{(\lambda iu\pm\mu (1-iu))\pi}{i}}\bigr)\Bigr)\\
            &=&\Bigl(\cos\bigl({\textstyle\frac{\lambda\pi}{j}}+(-\lambda\pm\mu )v\pi\bigr)\, ,\, 
                     \cos\bigl({\textstyle\frac{\mu\pi}{i}}+(\lambda\pm\mu )u\pi\bigr)\Bigr)\\
            &=&\Bigl(\cos\bigl({\textstyle\frac{\lambda\pi}{j}}\bigr)\, ,\, 
                     \cos\bigl({\textstyle\frac{\mu\pi}{i}}\bigr)\Bigr)\,\, .
\end{eqnarray*} 
This last equality uses the fact that $\lambda\equiv\mu\, ({\rm mod}\, 2)$. So part (c) is proved.

For part (d), reasoning as above gives
\[
k_1\equiv k_2\equiv\lambda\, ({\rm mod}\, j)\,\, ,\,\, k_1\equiv\mu\, ({\rm mod}\, i)\,\, ,\,\, 
{\rm and}\,\, k_2\equiv-\mu\, ({\rm mod}\, i)\,\, .
\]
Since $0<\lambda <j$, it follows that $k_1,k_2\not\in j\Z$. Similarly, $0<\mu <i$ implies that 
$k_1,k_2\not\in i\Z$. Therefore,
\[
\alpha^{-1}(S)\subseteq \Bigl\lbrace\,\,\cos \Bigl(\frac{k\pi}{ij}\Bigr)\,\,\Bigl|
\,\, 1\le k<ij\, ,\, k\not\in i\Z\, ,\, j\not\in i\Z\,\,\Bigr\rbrace\,\, .
\]
The equality now follows by comparing cardinalities of these two sets. $\Box$

\begin{remark} {\rm The much-studied curves $y^n-x^m$ may be viewed as a degeneration of the curves $T_n(y)-T_m(x)$. The Lin-Zaidenberg Theorem implies that, when $\gcd (m,n)=1$, the curve $y^n-x^m$ has a unique embedding in the plane up to algebraic equivalence; see \cite{Abhyankar.Sathaye.96}. In the case $m$ and $n$ are distinct prime numbers, Abhyankar and Sathaye \cite{Abhyankar.Sathaye.96} generalized this as follows.
\begin{quote} (Uniqueness Theorem) {\it If $m$ and $n$ are distinct prime numbers, and if 
\[
f=y^n-x^m+\sum_{in+jm<mn}a_{ij}x^iy^j
\]
with $a_{ij}\in\kk$, then $f$ has only one place at infinity and $f$ has a unique plane embedding, up to algebraic equivalence.}
\end{quote}
Therefore, the curves $T_n(y)-T_m(x)=0$ admit only one planar embedding, up to algebraic equivalence, when $m$ and $n$ are distinct primes.}
\end{remark}

\begin{remark} {\rm In the 1954 article \cite{Fiedler.Granat.54} (in Czech), Fiedler and Granat studied certain rational curves with a maximum number of nodes, including the family defined by $(T_n,T_{n+1})$, $n\ge 1$.}
\end{remark}

\begin{figure}[ht] 
  \centering
  \includegraphics[width=5in]{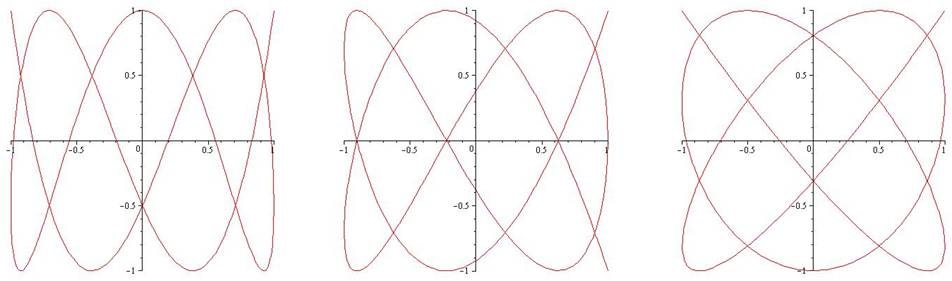}
  \caption{The curves $(T_3,T_8)$, $(T_4,T_7)$ and $(T_5,T_6)$}
  \label{fig:3-curves}
\end{figure}

\section{Line Embeddings $(T_i,T_j,T_k)$}

Given positive integers $i,j,k$, define the {\it pairwise greatest common divisor} of $i,j,k$ by
\[
{\rm pgcd}(i,j,k):=\max\{\gcd(i,j), \gcd (i,k), \gcd (j,k)\}\,\, .
\]

\subsection{Triples which Yield Embeddings}

\begin{proposition}\label{embedding} Let $\kk$ be a field of characteristice zero. Given integers $i,j,k\ge 2$, let $\phi$ denote the morphism $\phi =(T_i,T_j,T_k):\A^1\to\A^3$. Then $\phi$ is an embedding if and only if ${\rm pgcd}(i,j,k)=1$. 
\end{proposition}

\noindent {\it Proof.} Assume first that $\phi$ is an embedding. Set $d=\gcd (i,j)$, and write $i=ad$ and $j=bd$ for integers $a$ and $b$. Then
\[
T_1\in\kk [T_i,T_j,T_k]=\kk [T_a(T_d),T_b(T_d),T_k]\subset\kk [T_d,T_k]\,\, ,
\]
which implies that the morphism $(T_d,T_k):\A^1\to\A^2$ is an embedding. By the Epimorphism Theorem, either $k\vert d$ or $d\vert k$. 
If $k\vert d$, write $d=ck$. Then $T_1\in\kk [T_c(T_k), T_k]=\kk [T_k]$, which implies $k=1$, a contradiction. Therefore, $d\vert k$. 
Write $k=ed$. Then $T_1\in\kk [T_d, T_e(T_d)]=\kk [T_d]$, which implies $d=1$. In the same way, $\gcd (i,k)=\gcd (j,k)=1$. Therefore, ${\rm pgcd}(i,j,k)=1$.

Conversely, assume that ${\rm pgcd}(i,j,k)=1$. Then at most one of $i,j,k$ is even, and we may assume that $i$ is odd. Then it suffices to show the existence of positive integers $a,b,c$ such that $\vert aj-bk\vert = 1$ and $aj+bk=ci$. In this case, 
\[
2T_{aj}T_{bk}=T_{aj+bk}+T_{\vert aj-bk\vert}\quad\Rightarrow\quad T_1=2T_a(T_j)T_b(T_k)-T_c(T_i)\in\kk [T_i,T_j,T_k]\,\, ,
\]
and consequently $\phi$ is an embedding. 

It remains to show the existence of such integers $a$, $b$, and $c$. Since $\gcd (i,j)=1$, there exist positive integers $A$ and $B$ such that $\vert Aj-Bk\vert =1$. 
Thus, for any integer $x$, 
\[
\vert (A+kx)j-(B+jx)k\vert = 1\,\, .
\]
Consider the linear congruence $(A+kx)j+(B+jx)k\equiv 0\, ({\rm mod}\, i)$, or equivalently, 
\[
(2jk)x+(Aj+Bk)\equiv 0\, ({\rm mod}\, i)\,\, .
\]
Since $i$ is odd and ${\rm pgcd}(i,j,k)=1$, it follows that $\gcd (2jk,i)=1$, and thus $2jk\in\Z_i^*$. Choose $x$ such that $0\le x\le i-1$ and $x\equiv -(2jk)^{-1}(Aj+Bk)\, ({\rm mod}\, i)$. Then
the integers 
\[
a=kx+A\,\, ,\quad  b=jx+B\,\, ,\quad  {\rm and}\quad  c=(aj+bk)/i
\]
satisfy the required conditions. $\Box$ 

\begin{remark} {\rm Given integers $i,j,k$ with $\pgcd (i,j,k)=1$, the proof above gives an algorithm for finding integers $a,b,c$ such that $T_1=2T_a(T_j)T_b(T_k)-T_c(T_i)$.}
\end{remark}

\medskip

Note that if $1\in \{ i,j,k\}$, then $T_1\in\kk [T_i,T_j,T_k]$ and $(T_i,T_j,T_k)$ is an embedding. 
Combining this with {\it Prop.~\ref{embedding}} gives a complete description of which triples 
$(i,j,k)$ yield embeddings. We next investigate when two such embeddings are algebraically equivalent.

\subsection{Algebraically Equivalent Embeddings}

\noindent {\bf Example 1.} Over the field $\kk =\R$, the graph of the embedding $\phi:=(T_3,T_4,T_7)$ is given in {\it Fig.~\ref{fig:3,4,7}}. From this, it is evident that the image of $\phi$ is a trivial knot. That $\phi$ is a trivial embedding is seen as follows. Define $\epsilon\in GA_3(\R )$ by $\epsilon = (x,y,2xy-z):\R^3\to \R^3$. Since $T_1=2T_3T_4-T_7$, it follows that 
$\epsilon (T_3,T_4,T_7)=(T_3,T_4,T_1)$. The fact that $(T_3,T_4,T_1)$ is an algebraically trivial (rectifiable) embedding is a consequence of the following lemma. Note that rectifiability is {\it a priori} a much stronger property for an embedded real line than that of being topologically trivial.

\begin{figure}[ht] 
  \centering
  \includegraphics[width=2in]{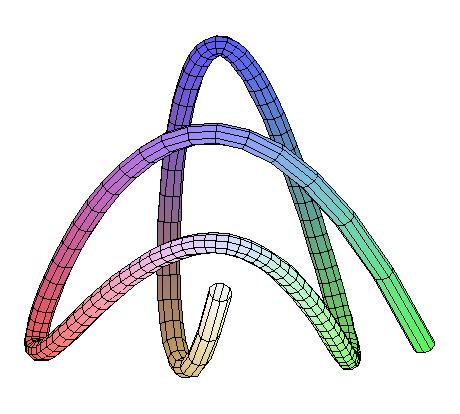}
  \caption{The unknot $(T_3,T_4,T_7)$}
  \label{fig:3,4,7}
\end{figure}

\begin{lemma}\label{trivial} If $i,j,k$ are positive integers and $1\in\{ i,j,k\}$, then $(T_i,T_j,T_k)$ is a trivial embedding.
\end{lemma}
\noindent {\it Proof.} Suppose $i=1$, and define $\beta\in GA_3(\kk )$ by $\beta = (x,T_j(x)-y,T_k(x)-z)$. Then 
\[
\beta (T_1,T_j,T_k)=(T_1,T_j(T_1)-T_j,T_k(T_1)-T_k)=(t,0,0)\,\,.
\]
$\Box$

\begin{definition} {\rm An {\it elementary involution} $\epsilon :\A^3\to\A^3$ is any polynomial map having one of the following forms:
\[
\begin{array}{l}
(x,z,y)\, ,\, (y,x,z)\, ,\, (z,y,x)\\
(f(y,z)-x,y,z)\, ,\, (x,f(x,z)-y,z)\, ,\, (x,y,f(x,y)-z)\\
(g(z)-x,h(z)-y,z)\, ,\, (g(y)-x,y,h(y)-z)\, ,\, (x,g(x)-y,h(x)-z)\,\, ,
\end{array}
\]
where $f\in\kk^{[2]}$ and $g,h\in\kk^{[1]}$. }
\end{definition}

\begin{definition} {\rm Define an equivalence relation on the set of ordered triples of positive integers by $(i,j,k)\sim (a,b,c)$ if and only if there exists a sequence
of elementary involutions $\epsilon_1,...,\epsilon_m$ $(m\ge 0)$ such that 
$(T_i,T_j,T_k)=\epsilon_1\epsilon_2\cdots\epsilon_m (T_a,T_b,T_c)$.}
\end{definition}

\begin{definition} {\rm Let $i,j$ be positive integers such that $\gcd (i,j)=1$. The {\it remnant} of $i$ and $j$ is
\[
\ll i,j\gg = \{ k\in\N\,\,\vert\,\, k>\max (i,j)\,\, ,\,\, k\not\in \langle i,j\rangle\,\, ,\,\, {\rm pgcd}(i,j,k)=1\}\,\, .
\]
Here, $\langle i,j\rangle$ denotes the semigroup of positive integers generated by $i$ and $j$. }
\end{definition} 
It is well known that when $\gcd (i,j)=1$, then the complement of $\langle i,j\rangle$ in the set of positive integers is finite, and its largest element is $ij-i-j$ (the Frobenius number of $i$ and $j$). Since ${\rm pgcd}(i,j)=1$, it follows that $\max\ll i,j\gg = ij-i-j$. 

\begin{definition} {\rm The ordered triple $(i,j,k)$ of positive integers is a {\it reduced triple} if and only if 
either {\bf (a)} $1\in\{ i,j,k\}$, or {\bf (b)} $2\le i<j$, $\gcd (i,j)=1$, and $k\in\ll i,j\gg$.}
\end{definition}

\begin{proposition}\label{reduced} Let $i,j,k$ be positive integers with ${\rm pgcd}(i,j,k)=1$. Then there exists a reduced triple $(I,J,K)$ such that $(i,j,k)\sim (I,J,K)$. In particular, the embeddings $(T_i,T_j,T_k)$ and $(T_I,T_J,T_K)$ are algebraically equivalent. 
\end{proposition}

\noindent {\it Proof.} We may assume $i<j<k$. The proposition obviously holds if $(i,j,k)$ is a reduced triple, so assume $(i,j,k)$ is non-reduced. Then $2\le i<j<k$ and $k\in\langle i,j\rangle$. 
Write $k=ai+bj$ for positive integers $a$ and $b$. Then 
\[
T_{\vert ai-bj\vert}=2T_{ai}T_{bj}-T_{ai+bj}=2T_a(T_i)T_b(T_j)-T_k\,\, .
\]
Define the elementary involution $\epsilon = (x,y,T_a(x)T_b(y)-z)$. Then 
$\epsilon (T_i,T_j,T_k)=(T_i,T_j,T_{\vert ai-bj\vert})$, which implies $(i,j,k)\sim (i,j,\vert ai-bj\vert )$.
In addition, ${\rm pgcd}(i,j,ai+bj)=1$ implies ${\rm pgcd}(i,j,\vert ai-bj\vert )=1$. Observe that the total degree of the triple has been reduced, i.e., $i+j+\vert ai-bj\vert < i+j+k$. Since the degree can only be reduced a finite number of times, we eventually arrive at a reduced triple $(I,J,K)$ equivalent to $(i,j,k)$. $\Box$

\begin{corollary}\label{2-case} If $(T_2,T_j,T_k)$ is an embedding, then it is a trivial embedding.
\end{corollary} 

\noindent {\it Proof.} By the preceding proposition, $(2,j,k)\sim (I,J,K)$ for a reduced triple $(I,J,K)$. If $1\in\{ I,J,K\}$, then the result follows from {\it Lemma~\ref{trivial}}. So assume $2\le I<J<K$. The proof of the proposition shows that $\min\{ 2,j,k\}\le\min\{ I,J,K\}$, which implies $I\le 2$. Therefore, $I=2$, and $J$ and $K$ must both be odd. But then $K\in\langle 2,J\rangle$, a contradiction. So this case cannot occur. $\Box$

\begin{corollary} If $(T_i,T_j,T_{i+j})$ is an embedding, then it is a trivial embedding. 
\end{corollary}

\noindent {\it Proof.} Assume $i<j$. In the proof of {\it Prop.~\ref{reduced}}, the first reduction replaces the triple $(i,j,i+j)$ with $(i,j,j-i)$, which is equivalent to $(i,j-i,j)$. 
Note that the triple $(i,j-i,j)$ maintains the form $(a,b,a+b)$. Therefore, this process produces a reduced triple $(I,J,I+J)$ equivalent to $(i,j,i+j)$. It follows that $1\in \{ I,J,I+J\}$.
The result now follows by {\it Lemma~\ref{trivial}}. $\Box$

\subsection{Knots Defined by Chebyshev Polynomials}

In this section, we assume $k=\R$. 
The results of the preceding section indicate that every embedding of the form $(T_i,T_j,T_k)$ is algebraically equivalent, via a sequence of elementary involutions, to one for which $k\in\ll i,j\gg$. 
Thus, for a fixed pair of integers $i$ and $j$ with $0<i<j$ and $\gcd (i,j)=1$, the set of reduced triples $(i,j,k)$ is finite. Moreover, {\it Lemma~\ref{trivial}} and {\it Cor.~\ref{2-case}} show that we need only consider cases in which $i\ge 3$. {\it Table 1} shows the remnant $\ll i,j\gg$ for relatively prime pairs $i,j$ such that $3\le i<j$ and $\frac{1}{2}(i-1)(j-1)\le 16$. Recall that 
$\frac{1}{2}(i-1)(j-1)$ is the number of nodes in the projection of $(T_i,T_j,T_k)$ onto the $xy$-plane, and that $\max\ll i,j\gg=ij-i-j$. In particular, 
$\ll i,j\gg\ne\emptyset$. 

\begin{center}
{\bf \footnotesize Table 1. Remnant for small degree pairs}\\
\begin{tabular}{|c|c|l|}
\hline
$(i,j)$ & {\footnotesize 1/2}$(i-1)(j-1)$ & $\qquad\quad \ll i,j\gg$\\ \hline\hline
(3,4) &3& $\{ 5\}$\\ \hline
(3,5) &4& $\{ 7\}$\\ \hline
(3,7) &6& $\{ 8,11\}$\\ \hline
(4,5) &6& $\{ 7,11\}$\\ \hline
(3,8) &7& $\{ 13\}$\\ \hline
(3,10) &9& $\{ 11,17\}$\\ \hline
(4,7) &9& $\{ 9,13,17\}$\\ \hline
(3,11) &10& $\{ 13,16,19\}$\\ \hline
(5,6)  &10& $\{ 7,13,19\}$\\ \hline
(3,13) &12& $\{ 14,17,20,23\}$\\ \hline
(4,9)  &12& $\{ 11,19,23\}$\\ \hline
(5,7)  &12& $\{ 8,9,11,13,16,18,23\}$\\ \hline
(3,14) &13& $\{ 19,25\}$\\ \hline
(5,8)  &14& $\{ 9,11,17,19,27\}$\\ \hline
(3,16) &15& $\{ 17,23,29\}$\\ \hline
(4,11) &15& $\{ 13,17,21,25,29\}$\\ \hline
(6,7)  &15& $\{ 11,17,23,29\}$\\ \hline
(3,17) &16& $\{ 19,22,25,28,31\}$\\ \hline
(5,9)  &16& $\{ 11,13,16,17,22,26,31\}$\\ \hline
\end{tabular}
\end{center}
{\it Table 1} shows that there are 63 reduced triples $(i,j,k)$ with $i\ge 3$ and $\frac{1}{2}(i-1)(j-1)\le 16$. {\it Table 2} gives details for the embeddings defined by triples up to $\frac{1}{2}(i-1)(j-1)\le 14$, in addition to a few other identified cases. 

\pagebreak

\begin{figure}[ht] 
  \centering
  \includegraphics[width=5in]{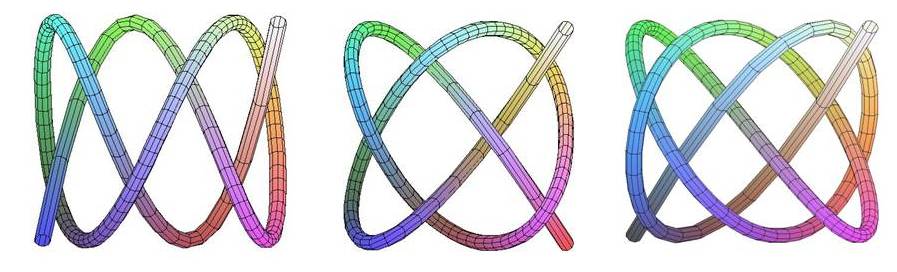}
  \caption{Knots for triples $(3,7,8)$, $(4,5,11)$ and $(5,7,9)$ }
  \label{fig:Cheb-group}
\end{figure}

\begin{center}
{\bf \footnotesize Table 2. Knot types for reduced triples of small degree}\\
\begin{tabular}{|c|c|c|c|c|c|c|c|c|}
\hline
   & $(i,j,k)$ & {\footnotesize $xy$-nodes} & {\footnotesize Knot type}& & & $(i,j,k)$ & {\footnotesize $xy$-nodes} & {\footnotesize Knot type} \\ \hline\hline
1. & (3,4,5) &3& $3_1$ & & 23. & (4,9,11) &12& $9(*)$ \\ \hline
2. & (3,5,7) &4& $4_1$ & & 24. & (4,9,19) &12& $10(*)$\\ \hline
3. & (3,7,8) &6& $5_1$ & & 25. & (4,9,23) &12& $12(*)$\\ \hline
4. & (3,7,11) &6& $6_3$ & & 26. & (5,7,8) &12& $4_1$\\ \hline
5. & (4,5,7) &6& $5_2$ & & 27. & (5,7,9) &12& $6_3$\\ \hline
6. & (4,5,11) &6& $6_2$ & & 28. & (5,7,11) &12&\\ \hline
7. & (3,8,13) &7& $7_7$ & &  29. & (5,7,13) &12&\\ \hline
8. & (3,10,11) &9& $7_1$ & & 30. & (5,7,16) &12& $9(*)$\\ \hline
9. & (3,10,17) &9& $9_{31}$ & &   31. & (5,7,18) &12&\\ \hline
10. & (4,7,9) &9& $7_5$ & &  32. & (5,7,23) &12& $12(*)$\\ \hline
11. & (4,7,13) &9& $8_7$ & &  33. & (3,14,19) & 13 & $11(*)$\\ \hline
12. & (4,7,17) &9& $9_{20}$ & &  34. & (3,14,25) & 13 & $13(*)$\\ \hline
13. & (3,11,13) &10& $8_3$ & & 35. & (5,8,9) & 14 & $7_3$ \\ \hline
14. & (3,11,16) &10& $9(*)$ & & 36. & (5,8,11) & 14 & $7_7$\\ \hline
15. & (3,11,19) &10& $10(*)$ & & 37. & (5,8,17) & 14 & $7_4$ \\ \hline
16. & (5,6,7) &10& $5_2$ & &  38. & (5,8,19) & 14 & $11(*)$\\ \hline
17. & (5,6,13) &10&  & & 39. & (5,8,27) & 14 & $14(*)$\\ \hline
18. & (5,6,19) &10& $10_{116}$ & & 40. & (3,16,17) & 15 & $11_1$\\ \hline
19. & (3,13,14) &12& $9_1$ & &  41. & $(6,7,11)$ & 15 & $8_{15}$\\ \hline
20. & (3,13,17) &12& $10(*)$ & & 42. & $(5,9,11)$ & 16 & $8_{12}$ \\ \hline
21. & (3,13,20) &12& $11(*)$& &  43. & $(7,8,9)$ & 21 & $7_5$\\ \hline
22. & (3,13,23) &12& $12(*)$& &  44. & $(9,10,11)$ & 36 & $9_{18}$\\ \hline

\end{tabular}
\end{center}
The method for determining knot types in {\it Table 2} was empirical, using a computer algebra system to construct the knot, then using Reidemeister moves and/or Jones polynomials to specify the type. It is known that if every singularity in the projection of a knot $K$ onto a plane is a node, and if the projection is reduced and alternating -- i.e., the crossings alternate between over- and under-crossings as the knot is traced out -- then the crossing number of $K$ equals the number of nodes in the projection; see \cite{Adams.01}, p.162. Thus, in some cases, it was possible to determine the crossing number of a knot in {\it Table 2}, even though the particular knot could not be identified. These cases are indicated by $(*)$ in the table. 
{\it Fig.~\ref{fig:Cheb-group}} illustrates three knots from {\it Table 2}.

\section{Nodal Curves and Crossing Sequences}

In this section, the assumption that the ground field is $\kk =\R$ continues.

\subsection{The Parity Property}
\begin{definition} {\rm A curve $C\subset\A^2$ is a {\it nodal curve} if it admits a parametrization $\alpha :\A^1\to\A^2$ such that $\alpha^{\prime}$ is non-vanishing and every singularity of $C$ is a node.}
\end{definition}

\begin{proposition}\label{parity} {\rm (Parity Property for Nodal Curves)} Let $C\subset\A^2$ be a nodal curve parametrized by $\alpha (t)$, and let $S\subset C$ be the set of nodes of $C$, where $\vert S\vert =n$. Suppose $\alpha^{-1}(S)=\{ t_1,...,t_{2n}\}$, where $t_1<t_2<\cdots <t_{2n}$. 
If $\alpha^{-1}(P)=\{ t_a,t_b\}$ for $P\in S$, then $a$ and $b$ are of opposite parity. 
\end{proposition}

\noindent {\it Proof.} If $b=a+1$, the result is obviously true, so assume $b-a\ge 2$. Let $K=\alpha ([t_a,t_b])\subset\R^2$, noting that since $\alpha (t_a)=\alpha (t_b)$, the restriction of $\alpha$ to $[t_a,t_b]$ factors as $[t_a,t_b]\to S^1\to K$, where $S^1$ denotes the unit circle. Let $\pi :H\to\R^2$ be the blow-up of $\R^2$ at the nodes of $K$, and let $C',K'\subset H$ be the proper transforms of $C$ and $K$, respectively. Then $K'$ is homeomorphic to $S^1$. Since $H$ is a simply connected open surface, $K'$ has a well-defined interior and exterior in $H$. Denote by $L\subset C'$ the proper transform of $\alpha (\R-[t_a,t_b])$.

Let $E=\sum_{i=1}^n E_i$ be the exceptional divisor of $H$, where each $E_i$ is irreducible. Then each $E_i$ intersects $K'$ in either one point or two points. Assume $E\cap K'=\{ Q_1,...,Q_r\}\cup\{ U_1,V_1,...,U_s,V_s\}$, where $E_i\cap K'=\{ Q_i\}$ for $1\le i\le r$, and 
$E_i\cap K'=\{ U_{i-r},V_{i-r}\}$ for 
$r+1\le i\le s$. Note that $\{ Q_1,...Q_r\}$ is the set of distinct points of $K'$ which are nodes of $C'$. Thus, 
each $Q_i$ is a point at which $L$ intersects $K'$ transversally, crossing either from the exterior to the interior of $K'$ (a point of entry), or from the interior to the exterior of $K'$ (a point of exit). To each point of entry, there corresponds a unique point of exit. Therefore, $r$ is an even integer, which implies that $r+2s$ is even. Since $r+2s$ is equal to the number of integers lying strictly between $a$ and $b$, it follows that $a$ and $b$ are of opposite parity. $\Box$

\begin{definition} {\rm Let $C\subset\A^2$ be a nodal curve. 
A {\it crossing sequence} for $C$ is a finite sequence $a_n$ of elements of $\{ -1,1\}$ with the following property.
Suppose $C$ is parametrized by $\alpha (t)$, and that $S\subset\A^2$ is the set of nodes of $C$, where $\vert S\vert =N$, and where 
$\alpha^{-1}(S)=\{ t_1,t_2,...,t_{2N}\}$ for $t_1<t_2<\cdots <t_{2N}$. 
Then $a_ma_n=-1$ whenever $\alpha (t_m)=\alpha (t_n)$. }
\end{definition}
A crossing sequence for the nodal curve $C$ encodes over- and under-crossing data for a knot $K$ having $C$ as regular projection. In this way, a nodal curve together with a crossing sequence determines a knot, which is unique up to knot type. 

One of the main implications of the Parity Property is the following.
\begin{corollary}
Let $C\subset\A^2$ be a nodal curve with $N$ nodes, and define the alternating sequence $a_n=(-1)^n$, $1\le n\le 2N$. 
Then $a_n$ is a crossing sequence for $C$.
Let $K$ be the knot determined by the alternating crossing sequence for $C$. Since $K$ admits an alternating regular projection onto $C$, its crossing number is $N$. 
\end{corollary}

\subsection{Alternating $(i,j)$-Knots}

\begin{definition} {\rm Given a pair of relatively prime positive integers $i<j$, the {\it alternating $(i,j)$-knot} is the knot $K$ determined by the alternating crossing sequence for the curve $(T_i,T_j)$. }
\end{definition}
Note that the crossing number of the alternating $(i,j)$-knot is $\frac{1}{2}(i-1)(j-1)$. This implies that the set of all alternating $(i,j)$-knots, taken over all possible pairs $(i,j)$, is t-infinite. {\it Table 3} below lists the alternating $(i,j)$-knots up to 10 crossings. 

\pagebreak

\begin{figure}[ht] 
  \centering
  \includegraphics[width=5in]{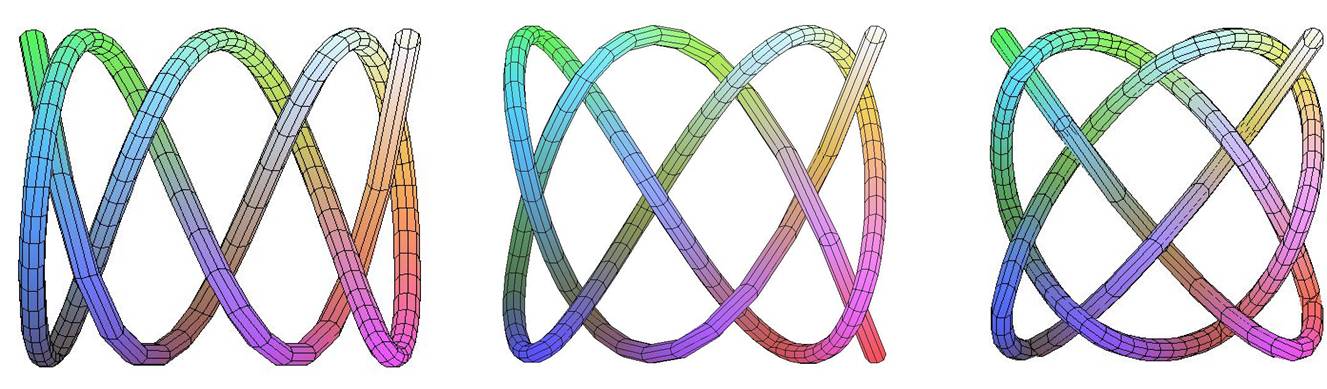}
  \caption{Alternating $(3,8)$, $(4,7)$, and $(5,6)$ knots}
  \label{fig:alt-group}
\end{figure}

\begin{center}
{\bf \footnotesize Table 3. Alternating $(i,j)$-Knots}\\
\begin{tabular}{|c|c|}
\hline
$(i,j)$ & {\footnotesize  Knot Type }\\ \hline\hline
(3,4) & $3_1$ \\ \hline
(3,5) & $4_1$ \\ \hline
(3,7) & $6_3$ \\ \hline
(4,5) & $6_2$ \\ \hline
(3,8) & $7_7$ \\ \hline
(3,10) & $9_{31}$ \\ \hline
(4,7) & $9_{20}$ \\ \hline
(3,11) & $10(*)$ \\ \hline
(5,6)  & $ 10_{116}$ \\ \hline
\end{tabular}
\end{center}
Three knots from {\it Table 3} are illustrated in {\it Fig.~\ref{fig:alt-group}}. 

Given coprime $i$ and $j$,
it is natural to look for a polynomial $h(t)$ such that $(T_i,T_j,h(t))$ parametrizes the alternating $(i,j)$-knot. 
In fact, $h(t)$ can be any polynomial which alternates in sign at the nodes of $(T_i,T_j)$; a necessary condition is that $\deg h(t)\ge ij-i-j$. However, we seek a closed form for $h(t)$ as a function of $i$ and $j$. We construct such $h(t)$ explicitly using Chebyshev polynomials of the second kind, $U_n(t)$.

As noted in {\it Lemma~\ref{roots}(c)}, the roots of $U_n(t)$ are 
\[
t=\cos \Bigl(\frac{k\pi}{n+1}\Bigr)\quad ,\quad k=1,...,n\,\, .
\]
Thus, given positive integers $m$ and $n$, $U_{m-1}$ and $U_{n-1}$ each divides $U_{mn-1}$.
{\it Lemma~\ref{identities}(h)} give the identity
\begin{equation}
U_{mn-1}=U_{m-1}(T_n)U_{n-1}=U_{n-1}(T_m)U_{m-1}\,\, .
\end{equation}
If $\gcd (m,n)=1$, then $\gcd (U_{m-1},U_{n-1})=1$, which implies that 
\begin{equation}
\frac{U_{mn-1}}{U_{m-1}U_{n-1}}
\end{equation}
is a polynomial of degree $(m-1)(n-1)$, and an even function, with roots
\[
t=\cos \Bigl(\frac{k\pi}{mn}\Bigr)\quad ,\quad k\in\Z_{mn}^*\,\, .
\]
In addition, each of these is a simple root, since roots of the polynomials $U_n$ are simple. 
The reader will note that the definition of the polynomials in (2) is quite similar to that of the cyclotomic polynomials $\Phi_{mn}$.  
\begin{proposition}\label{t-infinite} Let $m$ and $n$ be relatively prime positive integers, $3\le m<n$, and define 
\[
\phi = \Bigl(T_m,T_n,\frac{d}{dt}\frac{U_{mn-1}}{U_{m-1}U_{n-1}}\Bigr)\,\, .
\]
Then $\phi :\R^1\to\R^3$ is an embedding, and $\phi (\R^1)$ is the alternating $(m,n)$-knot.
\end{proposition}

\noindent {\it Proof.} Set 
\[
F(t)=\frac{U_{mn-1}}{U_{m-1}U_{n-1}}\,\, ,
\]
and let $t_1<t_2<\cdots <t_r$ be the roots of $F$, where $r=(m-1)(n-1)$. 
By {\it Prop.~\ref{curves}(d)}, the projection of $\phi (\R^1)$ onto the $xy$-plane has nodes at precisely the roots of $F(t)$. Since the roots of $F$ are simple, it follows that 
$F'(t_s)F'(t_{s+1})<0$ for each consecutive pair $t_s,t_{s+1}$ $(s=1,...,n-1)$. Suppose $\alpha (t_a)=\alpha (t_b)$ for $a<b$, where $\alpha =(T_m,T_n)$. 
Then by the Parity Property for Nodal Curves ({\it Prop.~\ref{parity}}), $a$ and $b$ are of opposite parity. It follows that $F'(t_a)F'(t_b)<0$. In particular, $F'(t_a)\ne F'(t_b)$. 
This implies that $\phi$ is one-to-one, since if $(x_0,y_0,z_0)$ is a point of self-intersection for the image of $\phi$, then $(x_0,y_0)$ is a node of $(T_m,T_n)$. Likewise, $\phi'$ is nowhere vanishing, since this is already true for $\alpha'$.
Therefore, $\phi$ is an embedding, and the projection of $\phi(\R^1)$ onto the $xy$-plane is an alternating projection. $\Box$

\medskip

\noindent {\bf Example.} Consider the case $m=3$ and $n=11$. Using the identity (1) above, and the fact that $U_2=4t^2-1$, we obtain
\[
F(t):=\frac{U_{32}}{U_2U_{10}}=\frac{U_2(T_{11})}{U_2}=\frac{4T_{11}^2-1}{4t^2-1}=\frac{2T_{22}+1}{4t^2-1}\,\, .
\]
The alternating $(3,11)$-knot is thus parametrized as follows. 
\begin{eqnarray*}
x=T_3&=& 4\, t^3-3\, t\\
y=T_{11}&=&1024\, t^{11}-2816\, t^9+2816\, t^7-1232\, t^5+220\, t^3-11\, t\\
z={\textstyle \frac{1}{64}}F'(t)&=&327680\, t^{19}-1548288\, t^{17}+3080192\, t^{15}-3347456\, t^{13}+2155008\, t^{11}\\
&&-832320\, t^9+185888\, t^7-21774\, t^5+1090\, t^3-15\, t
\end{eqnarray*}

\begin{remark} {\rm {\it Table 2} suggests that the knot $(T_i,T_j,T_k)$ is the alternating $(i,j)$-knot whenever $k=\max\ll i,j\gg = ij-i-j$, though we have not been successful in proving this. The difficulty in doing so is that, unlike $F'(t)$ above, the sign of $T_k$ does not alternate at the nodes of $(T_i,T_j)$.}
\end{remark}

\subsection{$(2,q)$ Torus Knots}

In this section, let $C_n\subset\A^2$ denote the nodal curve defined by $(T_3,T_{3n+1})$, $n\ge 1$. By {\it Prop.~\ref{curves}}, the number of nodes of $C_n$ is $3n$. Define a crossing sequence $a_m$, $1\le m\le 6n$, as follows. First,
\[
a_1=1\, ,\, a_2=-1\, ,\, a_3=-1\, , \quad {\rm and}\quad a_m=a_{m-1}a_{m-2}a_{m-3}\quad {\rm for}\quad 4\le m\le 2n\,\, .
\]
For the remaining elements, set
\[
a_{m+2n}=(-1)^{n+1}a_m \qquad {\rm and}\qquad a_{m+4n}=a_m\,\, \quad 1\le m\le 2n\,\, .
\]
The sign pattern for the first one-third of this sequence is either $+,-,-,+,+,-,-,...,-,-,+$ or $+,-,-,+,+,-,-,...,+,+,-$. If $n$ is odd, this pattern repeats twice, and if $n$ is even, the pattern alternates with its negation. 
\begin{proposition} For $n\ge 1$, the curve $C_n$ together with crossing sequence $a_m$, as defined above, define the $(2,2n+1)$ torus knot. \end{proposition}

\noindent {\it Proof.} Given a 2-tangle $V$, a pair of points $P,Q\in V$ is a {\it central pair} for $V$ if each belongs to distinct components of $V$ and lie between the crossings in a standard projection. For $n\ge 1$, let $K_{2n+1}$ be the $(2,2n+1)$ torus knot. For $n\ge 2$, $K_{2n+1}$ can be constructed from $K_{2n-1}$ as follows. Let $V\subset K_{2n-1}$ be a 2-tangle, and let $P,Q\in V$ be a central pair. Cut $K_{2n-1}$ at $P$ and $Q$, and insert an additional 2-tangle $W$. The resulting knot is $K_{2n+1}$. 

\begin{figure}[ht] 
  \centering
  \includegraphics[width=4in]{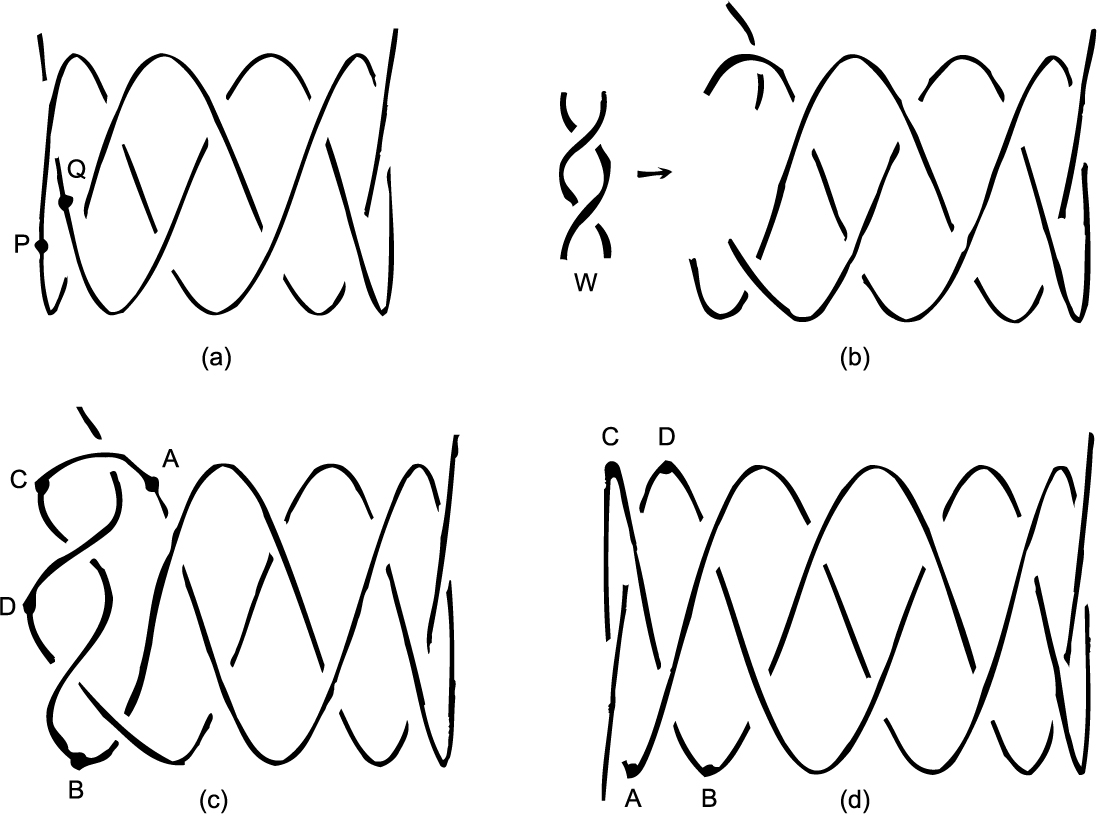}
  \caption{Splicing a 2-tangle into $C_n$}
  \label{fig:splicing}
\end{figure}

We now prove the proposition by induction on $n$, the case $n=1$ having been established in {\it Table 2}. For $n\ge 1$, assume $C_n$ is equivalent to $K_{2n+1}$. Let $V\subset C_n$ be the 2-tangle whose two crossings correspond to $a_{2n}=(-1)^n$ and $a_{2n+1}=(-1)^{n+1}$ in the crossing sequence $a_m$, and let $P,Q\in V$ be a central pair for $V$; see {\it Fig.~\ref{fig:splicing}(a)}. Cut $C_n$ at $P$ and $Q$, and insert a 2-tangle 
$W$, as per {\it Fig.~\ref{fig:splicing}(b)} and {\it (c)}. Now perform a twist as indicated in {\it Fig.~\ref{fig:splicing}(d)}, where $A,B,C,D$ represent corresponding points in {\it Fig.~\ref{fig:splicing}(c)} and {\it (d)}. The resulting knot is $C_{n+1}$. Since $C_{n+1}$ is obtained from $C_n$ by inserting a 2-tangle at the central pair of the preceding 2-tangle, beginning with the trefoil $C_1$, it follows by induction that, for all $n\ge 1$, $C_n$ is equivalent to $K_{2n+1}$. $\Box$

\begin{remark} {\rm {\it Table 2} suggests that $(T_3,T_{3n+1},T_{3n+2})$ parametrizes the $(2,2n+1)$ torus knot. This has been verified through $n=7$, though we have not succeeded in proving this for all $n$. The paper \cite{Koseleff.Pecker.pptb} shows that the $(2,2n+1)$ torus knot can be parametrized in degrees $(3,3n+1,3n+2)$, where the degree-3 component is $T_3$. }
\end{remark}

\section{Open Problems}

Various patterns can be observed in {\it Table 2}, suggesting several conjectures for the family of knots defined by $(T_i,T_j,T_k)$ for reduced triples $(i,j,k)$. 

\begin{conjecture} {\rm Every knot of the form $(T_i,T_j,T_k)$ for the reduced triple $(i,j,k)$, $i\ge 3$, is prime (hence non-trivial). }
\end{conjecture}

\begin{conjecture} {\rm Let $i,j$ be relatively prime positive integers, and let $k,\ell\in\ll i,j\gg$. If $k\ne\ell$, then $(T_i,T_j,T_k)$ and $(T_i,T_j,T_{\ell})$ define distinct knots.}
\end{conjecture} 


\begin{conjecture} 
{\rm Given odd $i\ge 3$, $(T_i,T_{i+1}, T_{i+2})$ parametrizes a knot with crossing number $i$.}
\end{conjecture}
\begin{conjecture} 
{\rm Given odd $i\ge 3$, $(T_i,T_{i+2}, T_{i+4})$ parametrizes a knot with crossing number $i+1$.}
\end{conjecture}


\noindent In addition:

\paragraph{Problem.} Given a reduced triple $(i,j,k)$, find an explicit formula for an invariant of $(T_i,T_j,T_k)$ as a function of $i,j,k$, for example, the crossing number. 

\medskip

\noindent Finally, we note that the standard three-term skein relations used to define invariant polynomials are similar to the recursive relation used to define Chebyshev polynomials. In fact, the Conway polynomial for the $(2,n)$ torus knot or link is precisely the (monic) Chebyshev polynomial of the first kind of degree $n-1$. It would be of interest to know if there is a connection between the invariant polynomials of a knot and the equations which define its parametrization.


\bibliographystyle{amsplain}

\bigskip

\noindent Department of Mathematics\\
Western Michigan University\\
Kalamazoo, Michigan 49008 USA\\
{\tt gene.freudenburg@wmich.edu}

\bigskip

\noindent Kalamazoo Area Math and Science Center\\
600 W Vine Street\\
Kalamazoo, Michigan 49008 USA\\
{\tt jennakay{\_}90@yahoo.com}

\end{document}